\documentstyle[12pt]{article}
 
\begin{document}


\def\a{\alpha}
\def\b{\beta}
\def\d{\delta}
\def\e{\epsilon}
\def\g{\gamma}
\def\k{\kappa}
\def\l{\lambda}
\def\o{\omega}
\def\t{\theta}
\def\s{\sigma}
\def\D{\Delta}
\def\L{\Lambda}

\def\G{{\cal G}}
\def\Gk{{\cal G}^{(2)}}
\def\R{{\cal R}}
\def\hR{{\hat{\cal R}}}
\def\C{{\bf C}}
\def\P{{\bf P}}

\def\uqgh{{U_q[gl(m|n)^{(1)}]}}
\def\uuqoh{{U_q[osp(1|2)^{(1)}]}}
\def\tuqoh{{U_q[osp(2|2)^{(2)}]}}
\def\tuqa22{{U_q[A^{(2)}_2]}}


\def\beq{\begin{equation}}
\def\eeq{\end{equation}}
\def\bea{\begin{eqnarray}}
\def\eea{\end{eqnarray}}
\def\ba{\begin{array}}
\def\ea{\end{array}}
\def\no{\nonumber}
\def\lt{\left}
\def\rt{\right}
\newcommand{\bq}{\begin{quote}}
\newcommand{\eq}{\end{quote}}
\newtheorem{Theorem}{Theorem}
\newtheorem{Definition}{Definition}
\newtheorem{Proposition}{Proposition}
\newtheorem{Lemma}[Theorem]{Lemma}
\newtheorem{Corollary}[Theorem]{Corollary}
\newcommand{\proof}[1]{{\bf Proof. }
        #1\begin{flushright}$\Box$\end{flushright}}

\newcommand{\sect}[1]{\setcounter{equation}{0}\section{#1}}
\renewcommand{\theequation}{\thesection.\arabic{equation}}


\title{\large \bf Drinfeld Basis And a Nonclassical Free Boson
    Representation of Twisted Quantum Affine Superalgebra $\tuqoh$} 
\author{{\large Wen-Li Yang~$^1$ and Yao-Zhong Zhang~$^2$}\\
\\
{\small $^1$ Institute of Modern Physics, Northwest University, 
Xian 710069, China}\\ 
{\small $^2$  Department of Mathematics, University of Queensland,
Brisbane, Qld 4072, Australia}}
\date{}
\maketitle

\vspace{10pt}

\begin{abstract}
We derive from the super RS algebra the Drinfeld basis 
of the twisted quantum affine superalgebra $\tuqoh$ by means of
the Gauss decomposition technique.  We explicitly construct a
nonclassical level-one representation of $\tuqoh$ in terms of two
$q$-deformed free boson fields.
\vskip 2cm

\end{abstract}

\vspace{10pt}

\sect{Introduction\label{intro}}
The algebraic analysis method based on infinite dimensional
non-abelian symmetries 
has proved extremely successful in both formulating and solving lower
dimensional integrable systems \cite{Is88}\cite{Jm94}. The key elements
of this method are infinite-dimensional highest weight 
representations of 
quantum (super) algebras and vertex operators. 
Drinfeld bases of quantum (super) algebras are of great importance in  
constructing the  infinite-dimensional representations 
 and vertex operators.The Drinfeld bases of quantum affine bosonic 
algebras were given by Drinfeld \cite{Dri88}.
For the case of quantum affine superalgebras, the Drinfeld bases have 
been known only for $\uqgh$ \cite{Yam96,Zha97,Cai98} and $\uuqoh$ 
\cite{Gou98,Din99}. In this paper, we derive 
the Drinfeld basis of the twisted 
quantum affine superalgebra $\tuqoh$ by using the super version of 
the Reshetikhin-Semenov-Tian-Shansky (RS) algebra \cite{Res90}
the  Gauss decomposition techinque of Ding-Frenkel \cite{Din93}. 
Moreover, we explicitly construct a nonclassical level-one
representation of $\tuqoh$ by means of two $q$-deformed free boson
fields.

\section{Drinfeld Basis of $\tuqoh$}

The symmetric Cartan matrix of the twisted affine superalgebra
$osp(2|2)^{(2)}$ is $(a_{ij})$ with $a_{11}=a_{22}=1,~
a_{12}=a_{21}=-1$. Twisted quantum affine superalgebra $\tuqoh$ is
a $q$-analogue of the universal enveloping algebra of $osp(2|2)^{(2)}$
and is generated by the Chevalley generators $\{E_i,~F_i,~K^{\pm 1}_i |
i=0,1\}$. The ${\bf Z}_2$-grading of the Chevalley generators is
$[E_i]=[F_i]=1,~ [K_i]=0,~ i=0,1$. The defining relations are
\bea
&&K_iK_j=K_jK_i,\no\\
&&K_jE_iK^{-1}_j=q^{a_{ij}}E_i,~~~~
  K_jF_iK^{-1}_j=q^{-a_{ij}}F_i,\no\\
&&\{E_i,F_j\}=\d_{ij}\frac{K_i-K^{-1}_i}{q-q^{-1}}
\eea
plus $q$-Serre relations which we omit. Here $\{X,Y\}\equiv XY+YX$.

$\tuqoh$ is a ${\bf Z}_2$-graded
quasi-triangular Hopf algebra endowed with the coproduct $\D$,
counit $\e$ and antipode $S$ given by
\bea
&&\D(E_i)=E_i\otimes 1+K_i\otimes E_i,~~~~
  \D(F_i)=F_i\otimes K^{-1}_i+1\otimes F_i,\no\\
&&\D(K_i)=K_i\otimes K_i,~~~~\e(E_i)=\e(F_i)=0,~~~~\e(K_i)=1,\no\\
&&S(E_i)=-K^{-1}_iE_i,~~~~ S(F_i)=-F_i K_i,~~~~
  S(K_i)=K^{-1}_i,~~~i=0,1.
\eea
The antipode $S$ is a ${\bf Z}_2$-graded algebra anti-homomorphism,
i.e. for homogeneous elements $a,b\in \tuqoh$, we have
$S(ab)=(-1)^{[a][b]}S(b)S(a)$.

We now give our definition of $\tuqoh$ in terms of Drinfeld generators.
\begin{Definition}
$\tuqoh$ is an associative algebra with unit 1 and
the Drinfeld generators: $X^\pm(z)$ and $\psi^\pm(z)$, a central
element $c$ and a nonzero complex parameter $q$. $\psi^\pm(z)$
are invertible. The gradings of the generators are: $[X^\pm(z)]=1$ and
$[\psi^\pm(z)]=[c]=0$. The relations are given by
\bea
&&\psi^\pm(z)\psi^\pm(w)=\psi^\pm(w)\psi^\pm(z),\no\\
&&\psi^+(z)\psi^-(w)=\frac{(z_+q+w_-)(z_-+w_+q)}
    {(z_++w_-q)(z_-q+w_+)}\psi^-(w)\psi^+(z),\no\\
&&\psi^\pm(z)^{-1}X^+(w)\psi^\pm(z)=\frac{z_\pm+w q}{z_\pm q+w}
    X^+(w),\no\\
&&\psi^\pm(z)X^-(w)\psi^\pm(z)^{-1}=\frac{z_\mp+w q}{z_\mp q+w}
    X^-(w),\no\\
&&(z+wq^{\pm 1})X^\pm(z)X^\pm(w)+(zq^{\pm 1}+w)X^\pm(w)X^\pm(z)
   =0,\no\\
&&\{X^+(z), X^-(w)\}=\frac{1}{(q-q^{-1})zw}\lt[
     \d(\frac{w}{z}q^{c})\psi^+(w_+)
     -\d(\frac{w}{z}q^{-c})\psi^-(z_+)\rt].\label{current}
\eea
\end{Definition}

Expand the currents in the form
\bea
&&X^\pm(z)=\sum_{n\in{\bf Z}} X^\pm_n\;z^{-n-1},\no\\
&&\psi^\pm(z)=\sum_{n\in{\bf Z}}\psi^\pm_n\;z^{-n}
  =K^{\pm 1}\exp\lt(\pm(q-q^{-1})\sum_{n>0}H_{\pm n}z^{\mp n}\rt).
\eea
In terms of modes $\{H_n\,|\,n\in {\bf Z}-\{0\},~
X^\pm_n\,|\,n\in{\bf Z}\}$ and $K$, the defining relations of $\tuqoh$
are given by
\bea
&&c:~{\rm central~ element},\no\\
&&[K, H_n]=0,~~~~KX^\pm_n K^{-1}=q^{\pm 1}X^\pm_n,\no\\
&&[H_n, H_m]=\d_{n+m,0}(-1)^n\frac{[n]_q[nc]_q}{n},~~~n\neq 0,\no\\
&&[H_n,X^\pm_m]=\pm\frac{1}{n}(-1)^n[n]_q
  q^{\mp \frac{c}{2}|n|}X^\pm_{n+m},~~~n\neq 0,\no\\
&&X^\pm_{n+1}X^\pm_m+q^{\pm 1}\,X^\pm_mX^\pm_{n+1}
  +q^{\pm 1}\,X^\pm_nX^\pm_{m+1}+X^\pm_{m+1}X^\pm_n=0,\no\\
&&\{X^+_n, X^-_m\}=\frac{1}{q-q^{-1}}\lt(q^{\frac{c}{2}(n-m)}
   \psi^+_{n+m}-q^{\frac{c}{2}(m-n)}\psi^-_{n+m}\rt).\label{mode}
\eea
Here and throughout $[i]_q=(q^i-q^{-i})/(q-q^{-1})$.

The Chevalley generators are obtained by the formulae:
\bea
&&K_1=K,~~~~E_1=X^+_0,~~~~F_1=X^-_0,\no\\
&&K_0=q^cK^{-1},~~~~E_0=X^-_1K^{-1},~~~~F_0=-KX^+_{-1}.\label{simple}
\eea
In terms of the Chevalley generators,
the Drinfeld generators can
be built up recursively by 
\bea
&&H_1=q^\frac{c}{2}K^{-1}(X^+_0X^-_1+X^-_1X^+_0),\no\\
&&H_{-1}=q^{-\frac{c}{2}}K(X^+_{-1}X^-_0+X^-_0X^+_{-1}),\no\\
&&X^\pm_{n+1}=\mp q^{\pm \frac{c}{2}}[H_1,X^\pm_n],~~~~
X^\pm_{-n-1}=\mp q^{\pm\frac{c}{2}}[H_{-1},X^\pm_{-n}],~~n\geq
  0
\eea
plus the formulae for $H_n,~H_{-n}~(n>0)$ given as follows
\beq
H_{\pm n}=\pm\frac{1}{q-q^{-1}}\sum_{p_1+2p_2+\cdots+np_n=n}
   \frac{(-1)^{\sum\,p_i-1}
       \;(\sum\,p_i-1)!}
   {p_1!\cdots p_n!}(K^{\mp 1}\psi^\pm_{\pm 1})^{p_1}
   \cdots (K^{\mp 1}\psi^\pm_{\pm n})^{p_n},\\
\eeq
where
\beq
\psi^\pm_{\pm n}=\pm(q-q^{-1})q^{\mp\frac{c}{2}(n-2)}
   \{X^\pm_{\pm n\mp 1},X^\mp_{\pm 1}\},~~n>0.
\eeq

\section{Derivation of Drinfeld Basis from Super RS Algebra}

Let us recall the definition of the super RS algebra. Let 
$R(z)\in End(V\otimes V)$, where $V$ is a ${\bf Z}_2$ graded vector space,
be a R-matrix which satisfies the 
graded Yang-Baxter equation 
\beq\label{rrr}
R_{12}(z)R_{13}(zw)R_{23}(w)=R_{23}(w)R_{13}(zw)R_{12}(z).
\eeq
We introduce \cite{Zha97,Gou98} 
\begin{Definition}\label{rs}: 
The super RS algebra $U(\R)$ is
generated by invertible L-operators $L^\pm(z)$, which obey the
relations
\bea
R({z\over w})L_1^\pm(z)L_2^\pm(w)&=&L_2^\pm(w)L_1^\pm(z)R({z\over
         w}),\no\\
R({z_+\over w_-})L_1^+(z)L_2^-(w)&=&L_2^-(w)L_1^+(z)R({z_-\over
         w_+}),\label{super-rs}
\eea
where $L_1^\pm(z)=L^\pm(z)\otimes 1$, $L_2^\pm(z)=1\otimes L^\pm(z)$
and $z_\pm=zq^{\pm {c\over 2}}$. For the first formula of
(\ref{super-rs}), the expansion direction of $R({z\over w})$ can be
chosen in $z\over w$ or $w\over z$, but for the second formula, the
expansion direction must only be in $z\over w$.
\end{Definition}
The multiplication rule for the tensor product is defined by
\beq
(a\otimes b)(a'\otimes b')=(-1)^{[b][a']}\,(aa'\otimes bb'),
\eeq
for homogeneous elements $a,~ b,~ a'$, $b'$ of $\tuqoh$.

In the following we apply the super RS algebra to derive the Drinfeld
basis of $\tuqoh$.
We take $R({z\over w})$ to be the 
R-matrix associated to the 3-dimensional representation $V$ of
$\tuqoh$. Let $v_1,~v_2,~v_3$ be the basis vectors of $V$ with
the ${\bf Z}_2$-grading $[v_1]=[v_3]=0$ and $[v_2]=1$.
It can be shown that the  R-matrix has the following form:
\beq
R({z\over w})=\left(
\begin{array}{ccccccccc}
1 & 0 & 0 & 0 & 0 & 0 & 0 & 0 & 0\\
0 & a & 0 & b & 0 & 0 & 0 & 0 & 0\\
0 & 0 & d & 0 & c & 0 & r & 0 & 0\\
0 & f & 0 & a & 0 & 0 & 0 & 0 & 0\\
0 & 0 & g & 0 & e & 0 & c & 0 & 0\\
0 & 0 & 0 & 0 & 0 & a & 0 & b & 0\\
0 & 0 & s & 0 & g & 0 & d & 0 & 0\\
0 & 0 & 0 & 0 & 0 & f & 0 & a & 0\\
0 & 0 & 0 & 0 & 0 & 0 & 0 & 0 & 1
\end{array}
\right), \label{r12}
\eeq
where
\bea
&&a=\frac{q(z-w)}{zq^2-w},~~~~b=\frac{w(q^2-1)}{zq^2-w},~~~~
  c=-\frac{q^{1/2}w(q^2-1)(z-w)}{(zq^2-w)(zq+w)},\no\\
&&d=\frac{q(z-w)(z+qw)}{(zq^2-w)(zq+w)},~~~~
  e=a-\frac{zw(q^2-1)(q+1)}{(zq^2-w)(zq+w)},\no\\
&&f=\frac{z(q^2-1)}{zq^2-w},~~~~
  g=-\frac{q^{1/2}z(q^2-1)(z-w)}{(zq^2-w)(zq+w)},\no\\
&&r=\frac{w^2(q-1)(q+1)^2}{(zq^2-w)(zq+w)},~~~~
  s=\frac{z^2(q-1)(q+1)^2}{(zq^2-w)(zq+w)}.
\eea

As in the non-super case \cite{Din93},
$L^\pm(z)$ allow the unique Gauss decomposition
\beq
L^\pm(z)=\left (
\begin{array}{ccc}
1 & 0 &  0\\
e^\pm_1(z) & 1 & 0\\
e^\pm_{3,1}(z) & e^\pm_2(z) & 1
\end{array}
\right )
\left (
\begin{array}{ccc}
k^\pm_1(z) & 0 & 0\\
0 & k^\pm_2(z) & 0\\
0 & 0 & k^\pm_3(z)
\end{array}
\right )
\left (
\begin{array}{ccc}
1 & f^\pm_1(z) & f^\pm_{1,3}(z) \\
0 & 1 &  f^\pm_2(z)\\
0 & 0 & 1
\end{array}
\right ),\label{l+-}
\eeq
where $e^\pm_{i,j}(z),~f^\pm_{j,i}(z)$ and $k^\pm_i(z) ~(i>j)$ are 
elements in the super RS algebra and $k^\pm_i(z)$ are invertible; and
$e_i^\pm(z)\equiv e^\pm_{i+1,i}(z),~f_i^\pm(z)\equiv f^\pm_{i,i+1}(z)$.
Let us define
\bea
X^+_i(z)&=&f^+_{i}(z_+)-f^-_{i}(z_-),\no\\
X^-_i(z)&=&e^-_{i}(z_+)-e^+_{i}(z_-).
\eea

By the definition of the super RS algebra and the Gauss decomposition
formula (\ref{l+-}), and after tedious
calculations parallel to those of the $U_q[osp(1|2)^{(1)}]$ case, 
we arrive at
\bea
k^\pm_i(z)k^\pm_j(w)&=&k^\pm_j(w)k^\pm_i(z),~~i,j=1,2,3,\no\\
k^+_1(z)k^-_1(w)&=&k^-_1(w)k^+_1(z),\no\\
k^+_3(z)k^-_3(w)&=&k^-_3(w)k^+_3(z),\no\\
\frac{z_\pm-w_\mp}{z_\pm q^2-w_\mp} k^\pm_1(z)k^\mp_2(w)&=&
    \frac{z_\mp-w_\pm}{z_\mp q^2-w_\pm}
    k^\mp_2(w)k^\pm_1(z),\no\\
\frac{(z_\mp-w_\pm)(z_\mp+w_\pm q)}{(z_\mp q^2-w_\pm)(z_\mp q+w_\pm)}
   k^\pm_1(z)k^\mp_3(w)^{-1}&=&
\frac{(z_\pm-w_\mp)(z_\pm+w_\mp q)}{(z_\pm q^2-w_\mp)(z_\pm q+w_\mp)}
   k^\mp_3(w)^{-1}k^\pm_1(z),\no\\
\frac{(z_\pm-w_\mp q^2)(z_\pm q+w_\mp)}{(z_\pm q^2-w_\mp)(z_\pm+w_\mp
   q)}k^\pm_2(z)k^\mp_2(w)&=&
  \frac{(z_\mp-w_\pm q^2)(z_\mp q+w_\pm)}{(z_\mp q^2-w_\pm)
  (z_\mp+w_\pm q)}k^\mp_2(w)k^\pm_2(z),\no\\
\frac{z_\pm-w_\mp}{z_\pm q^2-w_\mp} k^\pm_2(z)^{-1}k^\mp_3(w)^{-1}&=&
    \frac{z_\mp-w_\pm}{z_\mp q^2-w_\pm}
    k^\mp_3(w)^{-1}k^\pm_2(z)^{-1},\no
\eea
\bea
k^\pm_1(z)X^-_1(w)k^\pm_1(z)^{-1}&=&\frac{z_\pm q^2-w}{q(z_\pm-w)}
    X^-_1(w),\no\\
k^\pm_1(z)^{-1}X^+_1(w)k^\pm_1(z)&=&\frac{z_\mp q^2-w}{q(z_\mp-w)}
    X^+_1(w),\no\\
k^\pm_2(z)X^-_1(w)k^\pm_2(z)^{-1}&=&\frac{(z_\pm q^2-w)(z_\pm+wq)}
    {q(z_\pm-w)(z_\pm q+w)}X^-_1(w),\no\\
k^\pm_2(z)^{-1}X^+_1(w)k^\pm_2(z)&=&\frac{(z_\mp q^2-w)(z_\mp+wq)}
    {q(z_\mp-w)(z_\mp q+w)}X^+_1(w),\no\\
k^\pm_3(z)X^-_1(w)k^\pm_3(z)^{-1}&=&\frac{z_\pm+wq}{z_\pm q+w}
    X^-_1(w),\no\\
k^\pm_3(z)^{-1}X^+_1(w)k^\pm_3(z)&=&\frac{z_\mp+wq}{z_\mp q+w}
    X^+_1(w),\no\\
k^\pm_1(z)X^-_2(w)k^\pm_1(z)^{-1}&=&\frac{z_\pm q+w}{z_\pm+wq}
    X^-_2(w),\no\\
k^\pm_1(z)^{-1}X^+_2(w)k^\pm_1(z)&=&\frac{z_\mp q+w}{z_\mp+wq}
    X^+_2(w),\no\\
k^\pm_2(z)X^-_2(w)k^\pm_2(z)^{-1}&=&\frac{(z_\pm-wq^2)(z_\pm q+w)}
    {q(z_\pm-w)(z_\pm+wq)}X^-_2(w),\no\\
k^\pm_2(z)^{-1}X^+_2(w)k^\pm_2(z)&=&\frac{(z_\mp-wq^2)(z_\mp q+w)}
    {q(z_\mp-w)(z_\mp+wq)}X^+_2(w),\no\\
k^\pm_3(z)X^-_2(w)k^\pm_3(z)^{-1}&=&\frac{z_\pm-wq^2}{q(z_\pm-w)}
    X^-_2(w),\no\\
k^\pm_3(z)^{-1}X^+_2(w)k^\pm_3(z)&=&\frac{z_\mp-wq^2}{q(z_\mp-w)}
    X^+_2(w),\no
\eea
\bea
(z-wq^2)X^+_1(z)X^+_2(w)+q(z-w)X^+_2(w)X^+_1(z)
   &=&0,\no\\
q(z-w)X^-_1(z)X^-_2(w)+(z-wq^2)X^-_2(w)X^-_1(z)
   &=&0,\no\\
(z+wq^{\pm 1})X^\pm_1(z)X^\pm_1(w)+(zq^{\pm 1}+w)X^\pm_1(w)X^\pm_1(z)
   &=&0,\no\\
(z+wq^{\pm 1})X^\pm_2(z)X^\pm_2(w)+(zq^{\pm 1}+w)X^\pm_2(w)X^\pm_2(z)
   &=&0,\no
\eea
\bea
\{X^-_1(w), X^+_1(z)\}&=&(q-q^{-1})\lt[-\d(\frac{z}{w}q^c)k^+_2(z_+)
    k^+_1(z_+)^{-1}\rt.\no\\
& &~~~~~~~~\lt.+\d(\frac{z}{w}q^{-c})k^-_2(w_+)k^-_1(w_+)^{-1}\rt],\no\\
\{X^-_2(w), X^+_2(z)\}&=&(q-q^{-1})\lt[\d(\frac{z}{w}q^c)k^+_3(z_+)
    k^+_2(z_+)^{-1}\rt.\no\\
& &~~~~~~~~\lt.-\d(\frac{z}{w}q^{-c})k^-_3(w_+)k^-_2(w_+)^{-1}\rt],\no\\
\{X^-_2(w), X^+_1(z)\}&=&(q-q^{-1})q^{-{1\over 2}}
    \lt[-\d(-\frac{z}{w}q^{c-1})k^+_2(z_+)k^+_1(z_+)^{-1}\rt.\no\\
& &~~~~~~~~\lt.+\d(-\frac{z}{w}q^{-c-1})k^-_3(w_+)k^-_2(w_+)^{-1}\rt],\no\\
\{X^-_1(w), X^+_2(z)\}&=&(q-q^{-1})q^{1\over 2}\lt[\d(-\frac{z}{w}q^{c+1})
    k^+_2(-z_+q)k^+_1(-z_+q)^{-1}\rt.\no\\
& &~~~~~~~~\lt.-\d(-\frac{z}{w}q^{-c+1})
    k^-_3(-w_+q^{-1})k^-_2(-w_+q^{-1})^{-1}\rt],
  \label{x1-x2}
\eea
where
\beq
\d(z)=\sum_{l\in {\bf Z}}\,z^l
\eeq
is a formal delta function which enjoys the following properties:
\beq
\d(\frac{z}{w})=\d(\frac{w}{z}),~~~~~~
  \d(\frac{z}{w})f(z)=\d(\frac{z}{w})f(w).
\eeq

Defining the algebraic homomorphism
\bea
&&X^\pm(z)=z(q-q^{-1})\lt[X^\pm_1(z)+X^\pm_2(-zq^{-1})\rt],\no\\
&&\psi^-(z)=(1+q^{-{1\over 2}}-q^{1\over 2})\phi_1(z)-\phi_2(-zq^{-1}),\no\\
&&\psi^+(z)=\psi_1(z)-(1+q^{{1\over 2}}-q^{-{1\over 2}})\psi_2(-zq^{-1}),
\eea
where $\phi_i(z)=k^+_{i+1}(z)k^+_i(z)^{-1},~~
\psi_i(z)=k^-_{i+1}(z)k^-_i(z)^{-1},~~i=1,2$,
then we obtain from  (\ref{x1-x2}), the current
commutation relations (\ref{current}).

\section{Level Zero Representation}

We consider the evaluation representation $V_x$ of $\tuqoh$, where
$V$ is the 3-dimensional graded vector space with basis vectors
$v_1,~v_2$ and $v_3$. Let $e_{ij}$ be the $3\times 3$ matrix satisfying
$e_{ij}v_k=\d_{jk}v_i$. In the homogeneous gradation, the Chevalley
generators of $\tuqoh$ are represented on $V_x$ by
\bea
E_1&=&e_{12}-e_{23},~~~~F_1=e_{21}+e_{32},~~~~K_1=q^{e_{11}-e_{33}},\no\\
E_0&=&xq^{-1}(-e_{21}+e_{32}),~~~~F_0=x^{-1}q(e_{12}+e_{23}),~~~~
  K_0=q^{-e_{11}+e_{33}}.
\eea
Let $V_x^{*S}$  denote the dual module of $V_x$, 
defined by $\pi_{V^{*}}(a)=\pi_{V}(S(a))^{st}$. On $V^{*S}_x$,
the Chevalley generators are given by
\bea
E_1&=&-(q^{-1}e_{21}+e_{32}),~~~~F_1=qe_{12}-e_{23},~~~~
   K_1=q^{-e_{11}+e_{33}},\no\\
E_0&=&-xq^{-1}(e_{12}+q^{-1}e_{23}),~~~~F_0=x^{-1}q(-e_{21}+qe_{32}),~~~~
  K_0=q^{e_{11}-e_{33}}.
\eea
The following proposition can be proved by induction.
\begin{Proposition}: The Drinfeld generators are represented on $V_x$ by
\bea
H_m&=&-x^m\frac{[m]_q}{m}\lt[-(-1)^mq^{-m}e_{11}+\lt(1-(-1)^mq^m\rt)e_{22}
    +e_{33}\rt],\no\\
X^+_m&=&-x^m\lt(-(-1)^me_{12}+q^{-m}e_{23}\rt),\no\\
X^-_m&=&x^m\lt((-1)^me_{21}+q^{-m}e_{32}\rt),~~~~~~
    K=q^{e_{11}-e_{33}},
\eea
and on $V^{*S}_x$ by
\bea
H_m&=&-(-1)^mx^m\frac{[m]_q}{m}q^{-m}\lt[e_{11}+\lt(1-(-1)^mq^{-m}\rt)e_{22}
    -(-1)^mq^{-m}e_{33}\rt],\no\\
X^+_m&=&-(-1)^mx^mq^{-m}\lt(q^{-m-1}e_{21}+(-1)^me_{32}\rt),\no\\
X^-_m&=&(-1)^mx^mq^{-m}\lt(q^{-m+1}e_{12}-(-1)^me_{23}\rt),~~~~~~
    K=q^{-e_{11}+e_{33}}.
\eea
\end{Proposition}

\section{Nonclassical Free Boson Realization at Level One}

We use the notation similar to that of \cite{Awa94,Kim97,Zha98}.
Let us introduce two sets of bosonic oscillators $\{a_n,\;c_n,\;Q_{a},\;
Q_{c}|n\in{\bf Z}\}$ which 
satisfy the commutation relations
\bea
&&[a_n, a_m]=\d_{n+m,0}(-1)^n\frac{[n]_q^2}{n},~~~~~
  [a_0, Q_{a}]=1,\no\\
&&[c_n, c_m]=\d_{n+m,0}\frac{[n]_q^2}{n},~~~~~
  [c_0, Q_{c}]=1.\label{oscilators}
\eea
The remaining commutation relations are zero. Introduce the $q$-deformed
free boson fields
\bea
&&a(z;\k)=Q_{a}+a_0\ln z
   -\sum_{n\neq 0}\frac{a_n}{[n]_q}q^{\k |n|}z^{-n},\no\\
&&c(z)=Q_{c}+c_0\ln z-\sum_{n\neq 0}\frac{c_n}{[n]_q}
   z^{-n}
\eea
and  set
\beq
a_\pm(z)=\pm(q-q^{-1})\sum_{n>0}a_{\pm n}z^{\mp n}\pm a_0\ln q.
\eeq
Then
\begin{Theorem}\label{free-boson}: 
The Drinfeld generators of $U_q[osp(2|2)^{(2)}]$  at level one 
have the following nonclassical realization by the free boson fields 
\bea
&&\psi^{\pm}(z)=e^{a_\pm(z)},\no\\
&&X^{\pm}(z)=:e^{\pm a(z;\mp\frac{1}{2})}\;Y^{\pm}(z):F^\pm,\label{boson}
\eea
where $F^+=q^{1/2}+q^{-1/2},~~F^-=1/(q-q^{-1})$ and
\bea
&&Y^+(z)=e^{c(q^{1/2}z)}+e^{-c(-q^{-1/2}z)},\no\\
&&Y^-(z)=e^{c(-q^{1/2}z)}+
  e^{-c(q^{-1/2}z)}.
\eea
\end{Theorem}
This free boson realization is nonclassical in the sense that $X^-(z)$
defined by (\ref{boson}) does not have classical or $q\rightarrow 1$
limit.

\noindent{\it Proof.} We prove this theorem by checking that the
bosonized currents (\ref{boson})
satisfy the defining relations (\ref{current}) of $U_q[osp(2|2)^{(2)}]$ 
with $c=1$.
It is easily seen that the first two relations in (\ref{current})
are true by construction.
The third and fourth ones follow from the definition
of $X^\pm(z)$ and the commutativity between $a_n$ and $c_n$. So
we only need to check the last two relations in (\ref{current}).

We write
\beq
Z^\pm(z)=:e^{\pm a(z;\mp\frac{1}{2})}: .
\eeq
We obtain the operator products
\bea
Z^\pm(z)Z^\pm(w)&=&
  (z+q^{\mp 1}w)\;:Z^\pm(z)Z^\pm(w):,\no\\
Z^+(z)Z^-(w)&=&
  (z+w)^{-1}\;:Z^+(z)Z^-(w):.\no\\
Y^\pm(z)Y^\pm(w)&=&\pm q^{1/2} (z-w):e^{c(\pm q^{1/2}z)}e^{c(\pm q^{1/2}w)}:
  \mp q^{-1/2} (z-w):e^{- c(\mp q^{-1/2}z)}e^{- c(\mp q^{-1/2}w)}:\no\\
& & \pm \frac{q^{-1/2}}{z+q^{-1}w}:e^{c(\pm q^{1/2}z)}e^{- c(\mp q^{-1/2}w)}:
  \mp \frac{q^{1/2}}{z+qw}:e^{-c(\mp q^{-1/2}z)}e^{c(\pm q^{1/2}w)}:,\no\\
Y^+(z)Y^-(w)&=&q^{1/2} (z+w):e^{c(q^{1/2}z)}e^{c(-q^{1/2}w)}:
  -q^{-1/2} (z+w):e^{-c(-q^{-1/2}z)}e^{-c(q^{-1/2}w)}:\no\\
& &  +\frac{q^{-1/2}}{z-q^{-1}w}:e^{c(q^{1/2}z)}e^{-c(q^{-1/2}w)}:
  -\frac{q^{1/2}}{z-qw}:e^{-c(-q^{-1/2}z)}e^{c(-q^{1/2}w)}:.\no
\eea
Then the second last relation in (\ref{current}) is easily seen to be
true, and as to the last relation we have
\bea
\{X^+(z), X^-(w)\}&=&\frac{q^{1/2}+q^{-1/2}}{q-q^{-1}}:Z^+(z)Z^-(w):\no\\
& &\lt[\lt(\frac{q^{-1/2}}{(z+w)(z-q^{-1}w)}+\frac{q^{1/2}}{(w+z)
   (w-qz)}\rt):e^{c(q^{1/2}z)}e^{-c(q^{-1/2}w)}:\rt.\no\\
& &\lt.-\lt(\frac{q^{1/2}}{(z+w)(z-qw)}+\frac{q^{-1/2}}{(w+z)
   (w-q^{-1}z)}\rt):e^{-c(-q^{-1/2}z)}e^{c(-q^{1/2}w)}:\rt]\no\\
&=&\frac{1}{(q-q^{-1})zw}:Z^+(z)Z^-(w):\no\\
& &\lt[\lt(\d(-\frac{w}{z})-\d(\frac{w}{z}q^{-1})\rt)
  :e^{c(q^{1/2}z)}e^{-c(q^{-1/2}w)}:\rt.\no\\
& &\lt.-\lt(\d(-\frac{w}{z})-\d(\frac{w}{z}q)\rt)
  :e^{-c(-q^{-1/2}z)}e^{c(-q^{1/2}w)}:\rt]\no\\
&=&\frac{1}{(q-q^{-1})zw}
     \lt(\d(\frac{w}{z}q)
     -\d(\frac{w}{z}q^{-1})\rt):Z^+(z)Z^-(w):\no\\
&=&\frac{1}{(q-q^{-1})zw}
     \lt[\d(\frac{w}{z}q)\psi^+(wq^{1/2})
     -\d(\frac{w}{z}q^{-1})\psi^-(wq^{-1/2})\rt].
\eea
This completes the proof.

\section*{Acknowledgement}
Y.-Z.Z would like to thank Australia 
Research Council IREX programme for an Asia-Pacific Link 
Award and Institute of Modern Physics of Northwest University for 
hospitality.The financial support from Australian Research Council 
large, small and QEII fellowship grants is also gratefully 
acknowledged.


\end{document}